\documentclass{article}
\usepackage{amssymb,amsmath,color,enumerate}
\usepackage{wrapfig,epsfig,graphicx,curves}
\usepackage{makeidx}
 
\usepackage{latexsym}
\usepackage{pgf,tikz}
\usetikzlibrary{arrows}

\usepackage[applemac]{inputenc} % (pour les accents)
\renewcommand{\r}{\mathbb{R}}

\usepackage{graphicx}

\title{On the origin of Hilbert Geometry}
\author{Marc Troyanov - EPFL}

\begin{document}
\maketitle

\begin{abstract}
 In this brief essay we succinctly comment on the historical origin of Hilbert geometry. In particular, we give a summary of the letter in which David Hilbert  informs his friend and colleague Felix Klein about his discovery of this geometry.  The present paper is to appear in the Handbook of Hilbert geometry, (ed. A. Papadopoulos and M. Troyanov), European Mathematical Society, Z\"urich, 2014.

\end{abstract}
 
\vspace{0.7cm}

 Hilbert\index{Hilbert (David)} geometry, which is the main topic of this handbook, was born twelve decades ago, in the summer of 1894. Its birth certificate is a letter written by David Hilbert to Felix Klein, on August 14, 1894, sent from the village of  Kleinteich near  the city of Rauschen (today called Svetlogorsk), at the Baltic sea. The place was a fashionable vacation resort, situated some fourty kilometers north of Hilbert's hometown K\"onigsberg (now  Kaliningrad), where he was a Privatdozent\footnote{Hilbert was appointed professor at  G\"ottingen one year later.}. The mathematical part of this letter  was published one year later in the \emph{Mathematische Annalen}\footnote{The chief editor of this journal was Klein himelf.} under the title  \emph{Ueber die gerade Linie als k\" urzeste Verbindung zweier Punkte} (On the straight line  as shortest connection between two points)  \cite{Hilbert1895}.
 
\smallskip

Before discussing the content of this letter,  let us say a few words on the general situation of geometry in the $\mathrm{XIX}^\mathrm{{th}}$ century.  
At the dawn of the century,  and despite numerous criticisms concerning the axioms, Euclid's \emph{Elements} were still considered as  a model of mathematical rigor and an exposition of the ``true'' geometry.
The 1820s saw on one hand the creation of non-Euclidean geometry by Gauss, Bolyai and Lobachevsky, and
on the other hand, a revival and a rapid development of projective geometry by Poncelet, Gergonne, Steiner, M\"obius, Pl\"ucker and von Staudt. During this period, a  controversy developed among 
supporters of the synthetic versus the analytic methods and techniques in the subject.  In this context,  the \emph{analytic methods} are based on coordinates and algebraic relations for the description and the analysis of geometric figures, while the \emph{synthetic methods} are based on the incidence relations between points, lines, planes and other loci. 
The key to unify both sides of the controversy came from the work of von Staudt\index{von Staudt (Karl Georg Christian)} \cite{VonStaudt1847,VonStaudt1856}  who proved that numerical coordinates can be assigned to points on a line in a projective space through the use of purely synthetic methods and without using distances.  Von Staudt's   theorem implies that the cross ratio of four aligned points can be defined synthetically, and thus it is invariant under collineations (that is, transformations preserving lines and incidence).   

Recall that the cross ratio of four points $X,Y,Z,T$ is analytically defined as   
$$
  [X,Y,Z,T] =  \frac{(x-z)}{(x-t)} \cdot \frac{(y-t)}{(y-z)},
$$
where $x,y,z,t \in \r \cup \{\infty\}$ are  coordinates on the given line, corresponding to the points  $X,Y,Z,T$.
The synthetic definition of the cross ratio is more elaborate and it is based on an iteration of the construction of the harmonic conjugate of a point with respect to a  given pair of points (the harmonic conjugate corresponds to a cross ratio equal to $-1$ and can easily be constructed synthetically using complete quadrangles), together with ordering and density arguments. In fact, von Staudt was not able to fully prove his theorem since at his time the topological structure of the real numbers, in particular the completeness axiom, had yet to be clarified. This was done later on in Dedekind's \emph{Was sind und was sollen die Zahlen?} (1888).  A full exposition of von Staudt's construction is given in chapter VII of Veblen and Young's book \cite{VeblenYoung},
see also \cite[chap. VII]{Robinson1940}. We refer to \cite{Nabonnan2008,Voelke2008} for additional comments on von Staudt's work. 

\smallskip 
 
By the middle of the $\mathrm{XIX}^\mathrm{{th}}$ century, it was seen as a natural problem to try to define metric notions from projective ones, thus
reversing the analytic way of seeing geometry. Von Staudt's theorem allowed the use of the cross ratio in this adventure. The first step in this direction is usually attributed to Laguerre,\index{Laguerre (Edmond)} who proved in 1853 that the Euclidean angle $\varphi$ between two lines $\ell  ,m$ through the origin in a plane is given by 
\begin{equation}\label{eq.Laguerre}
 \varphi = \left| \frac{1 }{2i} \log [X,Y,U,V]\right|,  \tag{1}
\end{equation}
where  $X  \in  \ell $, $Y \in m $ and $U ,V$ are points representing  the ``isotropic lines'' of the plane\footnote{In fact the formula in Laguerre's paper does not explicitly involve the cross ratio, the present formulation seems to be due to F. Klein, see e.g.  \cite[p. 158]{Klein1926}.}. 
These are lines which are ``orthogonal to themselves'', meaning that if $U=(u_1,u_2)$ and $V=(v_1,v_2)$,
then 
$$
 u_1^2 + u_2^2 =   v_1^2   + v_2^2 =  0,
$$
or, equivalently, $u_2= iu_1$ and $v_2 = -iv_1$. To explain  Laguerre's formula, one may assume $x_1=y_1=u_1=v_1=1$; 
then $x_2 = \tan (\alpha)$, $y_2 = \tan (\beta)$, $u_2 =i$ and $v_2 =-i$. We then have
$$
  [X,Y,U,V] = \frac{(\tan (\alpha)-i)}{(\tan (\alpha)+i)} \cdot \frac{(\tan (\beta)+i)}{(\tan (\beta)-i)}
 = \mathrm{e}^{2i(\alpha-\beta)}. 
$$

 Klein  in \cite{Klein1871}, inspired by a formula published by Cayley\footnote{Cayley's formula also did not involve the cross ratio; it is written in terms of homogenous coordinates.} in \cite{Cayley},
observed later on that  Equation  (\ref{eq.Laguerre}) also gives a projective definition of the distance
in elliptic geometry, that is, the metric on $\mathbb{RP}^n$ for which the standard projection $\mathbb{S}^n \to \mathbb{RP}^n$ is a local isometry. Here $X$ and $Y$ are arbitrary distinct points in $\mathbb{RP}^n$  and $U$ and $V$ are
points in $\mathbb{CP}^n$ which are aligned with $X$ and $Y$ and belong to the quadric
 \begin{equation} \label{eq.absolute}   
 \sum_{i=0}^n u_i^2 = 0.  \tag{2}
\end{equation}
Klein also considered similar  geometries where the quadric (\ref{eq.absolute}) is replaced by an arbitrary quadric, which, following Cayley's terminology, he called \emph{the absolute}. In 1871,  Klein realized that choosing as the absolute the standard cone $\sum_{i= 1}^n u_i^2 - u_0^2= 0$, the metric
restricted to the interior of the cone is a model of Lobatchevsky's geometry. See \cite{CP} for a historical discussion.

 \smallskip

In 1888, M. Pasch\index{Pasch (Moritz)} proposed a new axiomatic foundation of geometry based on the primitive notion of segment 
(or betweenness) rather than the projective notion of  lines. Pasch's geometry is often called  \emph{ordered geometry}\index{odered geometry} and it greatly influenced Hilbert. 

\medskip

Let us now return to Hilbert's letter. The letter starts with a discussion of how geometry should be founded on three ``elements'' (in the sense of primitive notions) named points, lines and planes (Hilbert's aim is a discussion of three-dimensional geometry). These elements  should satisfy three groups of axioms, which he only shortly discusses.\index{axioms!Hilbert}

I. The first groups are axioms of incidence\index{axioms!incidence} describing the mutual relations among points, lines and planes. These axioms state  that two distinct points determine a unique line, that three non aligned points determine a unique plane, that
a plane containing two distinct  points contains the full line determined by those points, and that two planes cannot meet at only one point. Furthermore,
every line contains at least two points, every plane contains at least three non aligned points and the space contains at least four non coplanar points.

II. The second groups of axioms deals with ordering points\index{axioms!order} on a line and they have been proposed by M. Pasch in \cite{Pasch}.  They state: between two distinct points $A$ and $B$ on a line, there exists at least a third point $C$,    and given three points on a line, one and only one of them lies between the two others.
Given two distinct points $A$ and $B$, there is another point $C$ on the same line such  that $B$ is between $A$ and $C$.
Given four points on a line, we can order them as $A_1,A_2,A_3,A_4$ is such a way that if $h<i<k$, then $A_i$ lies
between $A_h$ and $A_k$.  If a plane $\alpha$ contains a line $a$, then the plane is separated in two half-planes  such that 
$A$ and $B$ belong to the same half-plane if and only if there is no point of the line $a$ between $A$ and $B$.

III. The axiom of continuity.\index{axioms!continuity}  Let $A_i$ be an infinite sequence of points  on a line $a$. If there exists a point  $B$  
such that $A_i$ lies between $A_k$ and $B$ as soon as $k \leq i$, then there exists a point $C$ on the line $a$ such that 
 $A_i$ lies between $A_1$ and $C$ for any $i >1$  and any point $C'$ with the same property lies between $C $
 and $B$. (We recognize an axiom stating that any  decreasing  sequence on a line that is bounded below has
 a unique greatest lower bound).
 
 \smallskip
 
 Hilbert then claims that the von Staudt theory of harmonic conjugates can be fully developed on
 this axiomatic basis (in a ``similar way'' as in Lindemann's book, \emph{Vorlesungen \"uber Geometrie}).
 It therefore follows that coordinates can be introduced. More precisely, to each point, one can
 associate three  real numbers $x,y,z$ and   each plane corresponds to a linear relation among the coordinates.
 Furthermore, if one interprets  $x,y,z$ as rectangular coordinates in the usual Euclidean space, then the points
 in our initial space correspond to a certain convex\footnote{Hilbert, following Minkowski, uses the expression \emph{nirgend concaven K\"orper}, 
 which means \emph{nowhere concave}. In other words, the domain is convex, but not   necessarily strictly convex.} body in Euclidean space. Conversely, the points in  an arbitrary convex domain  represent points in our initial space: \emph{our initial space is build from the points in a convex region of the Euclidean space.}

 \smallskip 
 
The axiomatic framework set here by Hilbert is interesting in several aspects. The first group of axioms belongs to what is now called 
\emph{incidence geometry}. Observe that Hilbert formulates them for a three-dimensional space --  at that time it was unusual to discuss geometry in an arbitrary dimension $n$ -- but note also that this restriction reduces the complexity of the theory, because it is not so simple to define and discuss a general notion of dimension in   pure incidence geometry. Also  Desargues theorem (which plays a key role in the introduction of coordinates)  need not be added as an axiom since it can be proved in this framework. 

The second group of axioms describes \emph{ordered geometry}, see  \cite{Robinson1940} and
\cite{Prenowitz} for readable accounts of this topics and  \cite{Pambuccian2011}  for an impressive historical
compilation of the literature since Pasch's original work.  Notice that the notion of line (and hence plane) can be
defined on the sole base of the primitive notion of segment, but Hilbert found it convenient  to  have lines and planes
as primitive notions also. Again, this reduces the complexity of the theory (and furthermore it is in Euclid's spirit).

\medskip
 
 In a more contemporary language, what Hilbert is claiming here is that  \emph{an abstract space $M$ consisting of points, together with a structure made of lines, planes and the betweenness relation  (for points on a line) satisfying the said axioms is isomorphic to a convex region $\Omega$  in $\r^3$ with the usual notions of lines and planes and  the condition  that $C$ lies between 
$A$ and $B$  if and only if $C \in [A,B]$.}  It is not clear whether the axiomatics proposed in Hilbert's letter is sufficient to rigorously prove  
this result, and the question is also not addressed in his \emph{Grundlagen} \cite{Hilbert1899}.  However it is indeed the case that a synthetic axiomatic
characterization of convex domains in $\r^n$ can be based on ordered geometry. A theorem of this type is given in W.A. Coppel's book, who calls it \emph{the fundamental theorem of ordered geometry}, see \cite[p. 173]{Coppel}.

  \smallskip 
  
  At this point of his letter, 
  Hilbert stresses  that  arbitrary convex bodies also appear in Minkowski's work on number theory.
  He then proceeds to define a notion of length for the segment $AB$ in his general geometry, which he assumes
  is now realized as the interior of a convex region in $\r^3$. He follows Klein's construction and defines 
  the distance between two distinct  points $A$ and $B$ to be\footnote{the usual convention now is to divide
  this quantity by $2$.}
\begin{equation}\label{def_Hdistance}
    d(A,B) = \log [X,Y,B,A] =\log \left(\frac{\overline{YA}}{\overline{YB}} \cdot \frac{\overline{XB}}{\overline{XA}}\right)
    ,  \tag{3}
\end{equation}
  where $X,Y$ are the two points on the boundary of the convex domain aligned with $A$ and $B$ in the order $X,A,B,Y$. Observe that $d(A,B) >0$ since $ \frac{\overline{YA}}{\overline{AB}},  \frac{\overline{XB}}{\overline{XA}} > 1$
  (if $A=B$ it is understood that  $d(A,B) =0$).
  
\smallskip  
  
Hilbert observes that the distance (\ref{def_Hdistance}) depends on the given convex domain: \emph{if $X$ approaches $A$ and $Y$ approaches $B$,
then the distance  $d(A,B)$ increases}. In other words, a smaller convex domain gives rise to larger Hilbert distances.
Hilbert then proceeds to give a proof of the triangle inequality.

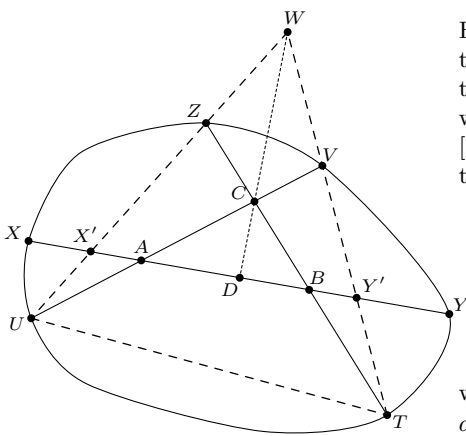
\begin{minipage}[r]{0.43\linewidth} \hspace{-2cm}
\begin{tikzpicture}[line cap=round,line join=round,>=triangle 45,x=0.7cm,y=0.7cm]
\clip(-2,-4.2) rectangle (10,5.5);
 \draw plot  [smooth cycle]  coordinates  { (0.24,-0.54)  
 (1.3,-1.79)  
 (4.64,-2.69) 
 (7,-2.38)  
 (8.18,-0.46)  
 (5.77,2.36)  
 (3.56,3.17)  
 (1.18,2.61)  
 (0.19,0.93)  }; 
\draw (0.19,0.93)-- (8.18,-0.46);
\draw (0.24,-0.54)-- (5.77,2.36); \draw (7,-2.38)-- (3.56,3.17);
\draw [line width=0.5pt,dash pattern=on 3pt off 3pt] (5.11,4.9)-- (0.24,-0.54);
\draw [line width=0.5pt,dash pattern=on 3pt off 3pt] (0.24,-0.54)-- (7,-2.38);
\draw [line width=0.5pt,dash pattern=on 3pt off 3pt] (7,-2.38)-- (5.11,4.9);
\draw [dash pattern=on 1pt off 1pt] (5.11,4.9)-- (4.2,0.23);
\begin{scriptsize}
\fill   (5.11,4.9) circle (1.5pt);
\draw (5.25,5.175) node {$W$};
\fill (0.24,-0.54) circle (1.5pt);
\draw  (-0.03,-0.64) node {$U$};
\fill   (7,-2.38) circle (1.5pt);
\draw  (7.25,-2.5) node {$T$};
\fill   (5.77,2.36) circle (1.5pt);
\draw  (5.93,2.56) node {$V$};
\fill   (3.56,3.17) circle (1.5pt);
\draw  (3.33,3.41) node {$Z$};
\fill   (8.18,-0.46) circle (1.5pt);
\draw  (8.41,-0.29) node {$Y$};
\fill   (0.19,0.93) circle (1.5pt);
\draw  (-0.1,1.08) node {$X$};
\fill   (4.2,0.23) circle (1.5pt);
\draw  (4.01,0.01) node {$D$};
\fill  (4.48,1.68) circle (1.5pt);
\draw  (4.19,1.84) node {$C$};
\fill   (2.33,0.56) circle (1.5pt);
\draw  (2.33,0.83) node {$A$};
\fill   (5.52,0) circle (1.5pt);
\draw  (5.67,0.24) node {$B$};
\fill   (1.37,0.73) circle (1.5pt);
\draw  (1.27,1.05) node {$X'$};
\fill   (6.42,-0.15) circle (1.5pt);
\draw  (6.75,0.1) node {$Y'$};
\end{scriptsize}
\end{tikzpicture}
\end{minipage}  
\begin{minipage}[l]{6cm} \vspace{-0.7cm} 
{\small Here is the argument with Hilbert's notation:  
From the invariance of the cross ratio with respect to the perspective at $W$,
we have  $ [U,V,C,A] = [X',Y',D,A]  $  and  $[Z,T,B,C] =  [X',Y',B,D]$.
Multiplying these identities gives
\begin{align*}
  [U,V,C,&A] \cdot  [Z,T,B,C] 
  \\ &=  [X',Y',D,A] \cdot [X',Y',B,D] 
  \\ &=  [X',Y',B,A] 
  \\ &\geq   [X ,Y,B,A], 
\end{align*}
which is equivalent to $d(A,C)+d(C,B) \geq d(A,B)$.
}
\end{minipage}

Hilbert also observes that the triangle inequality degenerates to an equality when the three points are aligned with $C$ between $A$ and $B$.
Furthermore, he discusses necessary and sufficient conditions for the equality in the triangle inequality. He shows that a non-degenerate triangle exists for which the sum of two sides is equal to the third if and only if there exists a plane which meets the boundary of
the convex domain on two segments not on the same line (these would be the segments $[V,T]$ and $[Z,U]$ on the figure). He quotes, as a particularly interesting example,  the case where the convex domain is bounded by a tetrahedron. 

As a final word, Hilbert points out that he always assumed the given convex body to be bounded, and this hypothesis implies that his
geometry does not satisfy the Euclidean parallel postulate.

\medskip

Let us conclude with two remarks. We first   mention that Hilbert came back at least  at  two occasions on the mathematical subjects discussed in his letter.
First in his book on the foundation of geometry \cite{Hilbert1899}, the book in which he systematically develops the axiomatics of Euclidean 3-space,
starting with the same axioms of incidence, order and continuity, to which he adds some axioms on congruence and the 
parallel postulate. The second occasion is the 1900 International Congress in Paris where he states the famous
Hilbert  problems. Problem IV concerns   {\it the construction and systematic treatment of the geometries} in
convex domains for which straight lines are the shortest, that is, the case where the triangle inequality is an equality in the
case of three aligned points.  Minkowski and Hilbert gometries are natural examples of such geometries.
We refer to the discussion in the last chapter of this Handbook \cite{Papadopoulos-Hilbert}.

\smallskip

The second and last remark is about how Hilbert understood the word \emph{geometry}. For
Hilbert, a geometry was not conceived as an abstract metric space\footnote{The notion of abstract metric space  first appeared
in the 1906  thesis of M. Fréchet \cite{Frechet}.} satisfying some specific axioms, but rather, as we saw, a geometry is a system made of points, lines and planes subject to a system of consistent interrelations
wich are accepted as axioms. A notion of distance between points in  a given geometry has then to be constructed from the given data and axioms and the properties of the distance, including the triangle inequality, is then a theorem that needs a proof rather than an axiom or a part of the initial definition. This is also in the spirit of Euclid's Element, where the triangle inequality is proved in Book 1, Proposition 20.
 
The subject of \emph{metric geometry} has been developed  since the 1920s by a long list of mathematicians including Hausdorff, Menger, Blumenthal, Urysohn, Birkhoff, Busemann,  Alexandrov, and  others. Since the work of  Gromov in the 1980s, metric geometry is seen as a topic of major importance in geometry.

\medskip

\textbf{Acknowledgement.}   The author thanks A. Papadopoulos and K.-D. Semmler for carefully reading the manuscript and providing useful comments.

\vfill

{\small Marc  Troyanov, Section de Math{\'e}matiques, \\ \'Ecole Polytechnique F{\'e}d\'erale de
Lausanne, \\ SMA--Station 8,  \\ CH-1015 Lausanne - Switzerland \\ marc.troyanov@epfl.ch}

\end{document}